# CARD SHUFFLING AND DIOPHANTINE APPROXIMATION

By Omer Angel ,[1] Yuval Peres[2] and David B. Wilson

*University of British Columbia and University of Toronto, University of California at Berkeley and Microsoft Research, and Microsoft Research*

The "overlapping-cycles shuffle" mixes a deck of $n$ cards by moving either the $n$th card or the $(n-k)$th card to the top of the deck, with probability half each. We determine the spectral gap for the location of a single card, which, as a function of $k$ and $n$, has surprising behavior. For example, suppose $k$ is the closest integer to $\alpha n$ for a fixed real $\alpha \in (0,1)$. Then for rational $\alpha$ the spectral gap is $\Theta(n^{-2})$, while for poorly approximable irrational numbers $\alpha$, such as the reciprocal of the golden ratio, the spectral gap is $\Theta(n^{-3/2})$.

## 1. Introduction.

The overlapping-cycles shuffle has two positive integer parameters $n$ and $k \in [1,n)$, and mixes a deck of $n$ cards by flipping a fair coin to decide whether to move the $n$th card or the $(n-k)$th card to the top of the deck. The case where $k=1$ was first studied by Hildebrand [6], who showed that the shuffle mixes within $O(n^3 \log n)$ time, and Wilson [8] who showed that it takes at least $(1-o(1))/(8\pi^2)n^3 \log n$ shuffles to mix. Jonasson [7] first studied the overlapping-cycles shuffle that we consider here, obtaining a lower bound of order $n^3/k^2 \log n$, but for this shuffle there was no matching upper bound, and indeed there cannot be one, since, as Jonasson showed, when $k = n/2$ the shuffle takes at least order $n^2$ time to mix.

We would like to understand the mixing time of this overlapping-cycle shuffle, but the behavior is sufficiently complicated that we settle for the more modest goal of understanding how a single card behaves under this shuffle, and in particular, what the spectral gap is for the Markov chain induced by a single card. The behavior of the spectral gap is a surprisingly complicated function of $n$ and $k$ (see Figures 1 and 2), which depends to a

Received July 2007; revised September 2007.
[1]Supported in part by a PIMS postdoctoral fellowship.
[2]Supported in part by NSF Grant DMS-06-05166.

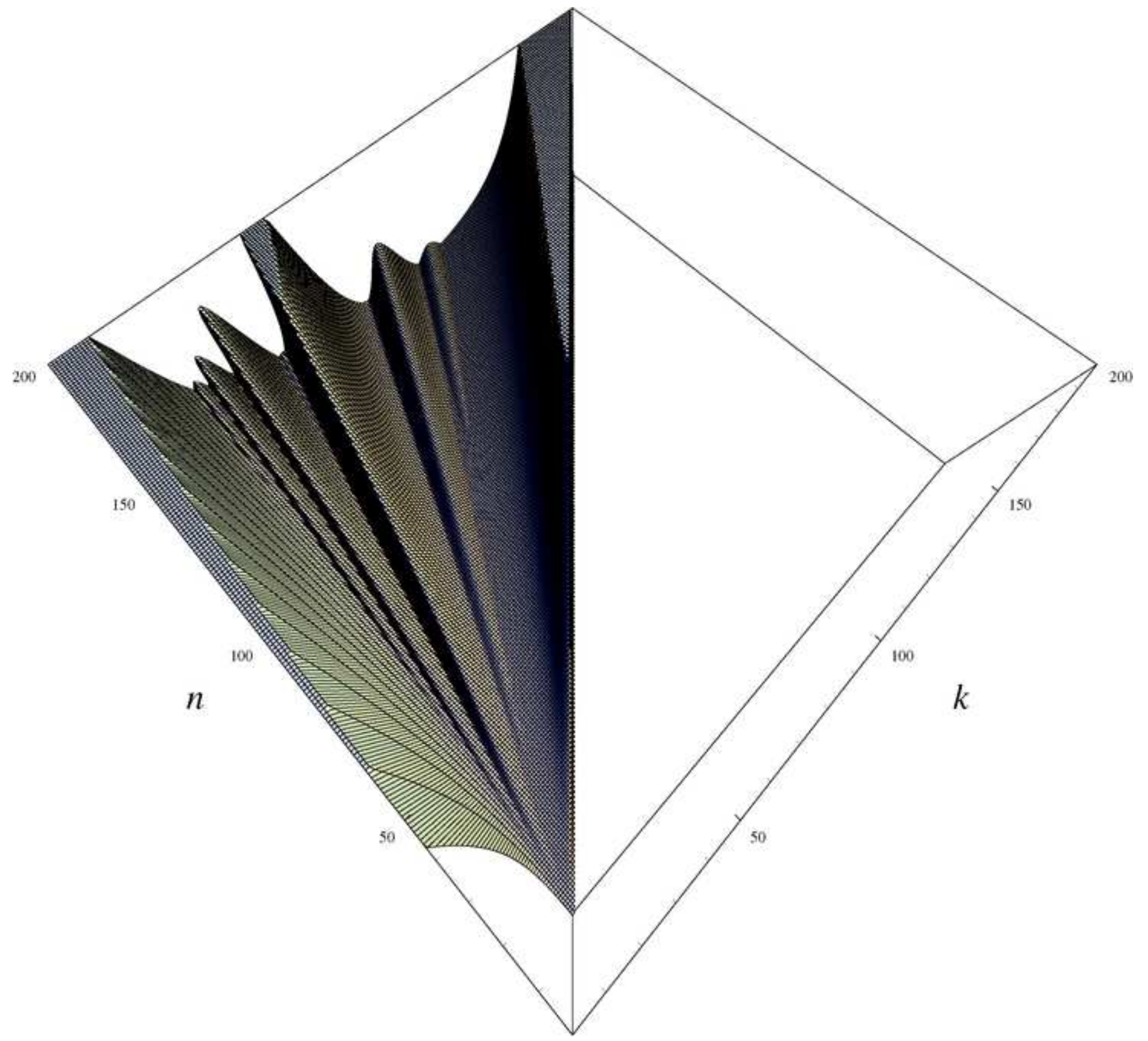

Fig. 1. *The relaxation time as a function of $n$ and $k$, for $1 \leq k < n \leq 200$.*

great extent on the Diophantine approximations of $k/n$. For example, suppose $\alpha$ is a fixed number between 0 and 1, and $k = \lfloor \alpha n \rceil$, that is, the closest integer to $\alpha n$. Then the inverse spectral gap (also known as the *relaxation time*) is of order $n^2$ when $\alpha$ is rational, but when $\alpha$ is a badly approximable irrational number, such as the golden ratio minus 1, the relaxation time is of order $n^{3/2}$.

1.1. *Heuristics.* Before giving more precise results below, we start with a heuristic discussion of why the relaxation time of a single card should fluctuate between order $n^2$ and order $n^{3/2}$ (when $k$ is of order $n$) according to the Diophantine approximations of $k/n$. We mention that Jonasson [7] already understood that the relaxation time is of order $n^2$ when $k/n$ is near a simple rational, but conjectured that the relaxation time is of order $n \log n$ for typical $k$.

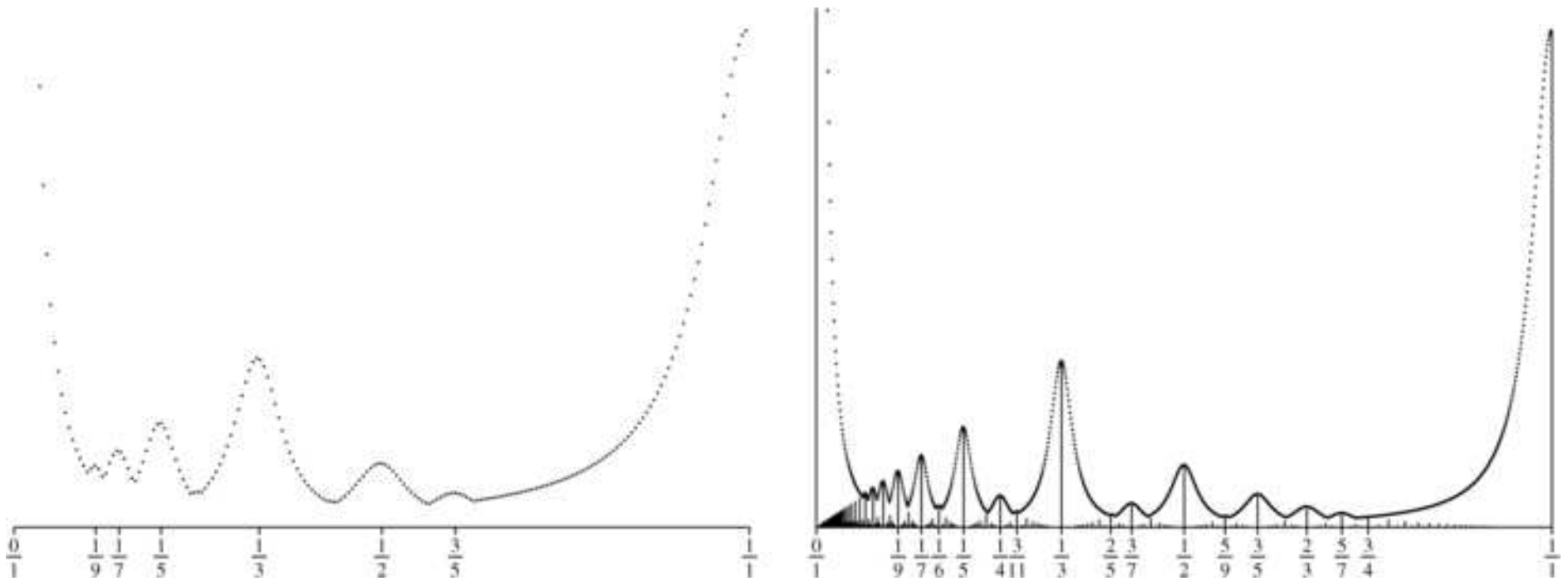

FIG. 2. *Shown here is the inverse spectral gap for the position of the ♠A for the $n,k$ shuffle; the left picture plots the inverse gap as a function of $k/n$ when $n=200$ and $k$ is varied from 1 to 199, while the right picture is for $n=1000$ with $k$ varying between 1 and 999. There are spikes when $k/n$ is well approximated by a simple rational number, and these spikes become more pronounced for larger $n$, with more of them visible. The predicted locations and heights of the peaks are shown in the $n=1000$ picture. At the tops of the spikes the inverse spectral gap is $\Theta(n^2)$, except for the spike at $\frac{0}{1}$ where it is $\Theta(n^3)$. At the low-lying areas between the spikes, the inverse spectral gap is $\Theta(n^{3/2})$.*

Suppose that the ♠A was initially at the top of the deck. Where will it be after $t$ steps? Denote by $s$ the number of visits of the ♠A to location 1 (*after* time 0), by $\widetilde{s}$ the number of times it jumped from location $n-k$ to 1, and by $\ell$ its final location, which we assume to be less than $n-k$ for simplicity. Clearly, $s=O(t/n)$ and $\widetilde{s}=s/2\pm O(\sqrt{s}\,)$. We get another approximate identity by noting that the card progresses at exactly unit speed from 1 to $n-k$ and then either jumps to 1 or moves at approximately speed 1/2 from $n-k$ to $n$. There are random fluctuations in how long the traversal from position $n-k$ to $n$ takes. We then have

$$t=\ell-1+(n-k)s+2k(s-\widetilde{s})\pm\Theta(\sqrt{2k(s-\widetilde{s})}\,).$$

But it is likely that $s-\widetilde{s}\approx t/(2n)$, and we assumed $k$ is of order $n$, so this error term is $\pm\Theta(\sqrt{t}\,)$. Thus, we may rewrite this in the form

$$\ell=k(2\widetilde{s}-s)-sn+t\pm\Theta(\sqrt{t}\,).$$

In other words, at time $t$, the location $\ell\in\{1,\dots,n\}$ of the ♠A is in one of

$$2\widetilde{s}-s=\pm\Theta(\sqrt{s}\,)=\pm\Theta(\sqrt{t/n}\,)$$

intervals $I_j$ of length $\Theta(\sqrt{t}\,)$ each; these intervals are centered at the first $\Theta(\sqrt{t/n}\,)$ multiples $\pm jk$ of $k$ modulo $n$. A necessary condition for mixing is that these intervals cover most of $\{1,\dots,n\}$, which requires $n=\Theta(t/\sqrt{n}\,)$, that is, $t$ must be at least of order $n^{3/2}$. If the multiples $\{jk:|j|\le\Theta(\sqrt{t/n}\,)\}$ of $k$ are close to uniformly distributed modulo $n$, and $t > Cn^{3/2}$, then the intervals $\{I_j : |j| \le \Theta(\sqrt{t/n})\}$ cover most of $\{1,\dots,n\}$. This happens when, for example, $n/k$ is close to the golden ratio $\frac{1+\sqrt{5}}{2}$, so in this case one might expect the relaxation time to be of order $n^{3/2}$. On the other hand, if $k/n$ is very close to a simple rational $p/q$, then there are only $q$ distinct intervals $I_j$ of length $\Theta(\sqrt{t})$ each; for these intervals to cover most of $\{1,\dots,n\}$, the time $t$ must be at least of order $(n/q)^2$.

1.2. *Statement of results.* We define $\operatorname{gap}(n,k)$ to be the spectral gap of a single card when the deck is shuffled according to the overlapping-cycles shuffle with parameters $n$ and $k$.

THEOREM 1. *Suppose that $p$ and $q$ are integers such that $0 < p/q < 1$ and $\gcd(p,q) = 1$. Then for fixed $p$ and $q$, when $n$ is large enough and $k = \lfloor \frac{p}{q} n \rceil$, the spectral gap is*

$$\operatorname{gap}(n,k) = (1 + o(1)) \frac{\pi^2 pq}{2n^2} \times \begin{cases} 1, & \text{if } p \equiv q \bmod 2, \\ 4, & \text{if } p \not\equiv q \bmod 2. \end{cases}$$

It is interesting to consider the spectral gap when $n$ is large but then held fixed while $k$ is varied, as in Figure 2. We prove the following bounds:

THEOREM 2.

(i) *For any $k \in [1,n)$, we have $\operatorname{gap}(n,k) \le (2\pi^2 + o(1)) n^{-3/2}$.*
(ii) *For any $k \in [1,n)$, we have $\operatorname{gap}(n,k) \ge (4\pi^2 + o(1)) n^{-3}$.*
(iii) *For almost all values of $k \in [1,n)$, we have that $\operatorname{gap}(n,k)$ is of order $n^{-3/2}$. More precisely, when $n$ is large enough, for any $\delta > 0$, for all but $\delta^{2/3} n$ values of $k \in [1,n)$, we have $\operatorname{gap}(n,k) \ge (\pi^2/16 + o(1)) \delta n^{-3/2}$.*

*Here the $o(1)$ terms tend to $0$ as $n \to \infty$ independent of $k$ or $\delta$.*

The lower bound is tight, but as discussed in the open problem section, the upper bound can almost certainly be improved by a factor of $\sqrt{3}$.

As the reader may have noticed already, the spikes in Figure 2 are actually shaped like narrow bells. These bells are shown more explicitly in Figure 3. Our next theorem characterizes the shape and width of these bells, by giving the spectral gap when $k$ is close to but not exactly $\lfloor (p/q) n \rceil$.

THEOREM 3. *Suppose that $p$ and $q$ are integers such that $0 < p/q \le 1$ and $\gcd(p,q) = 1$. Then for fixed $p$ and $q$, when $n$ is large enough and*

$$|k - (p/q)n| \le n^{3/4} \times \text{(constant depending on } p \text{ and } q\text{)},$$

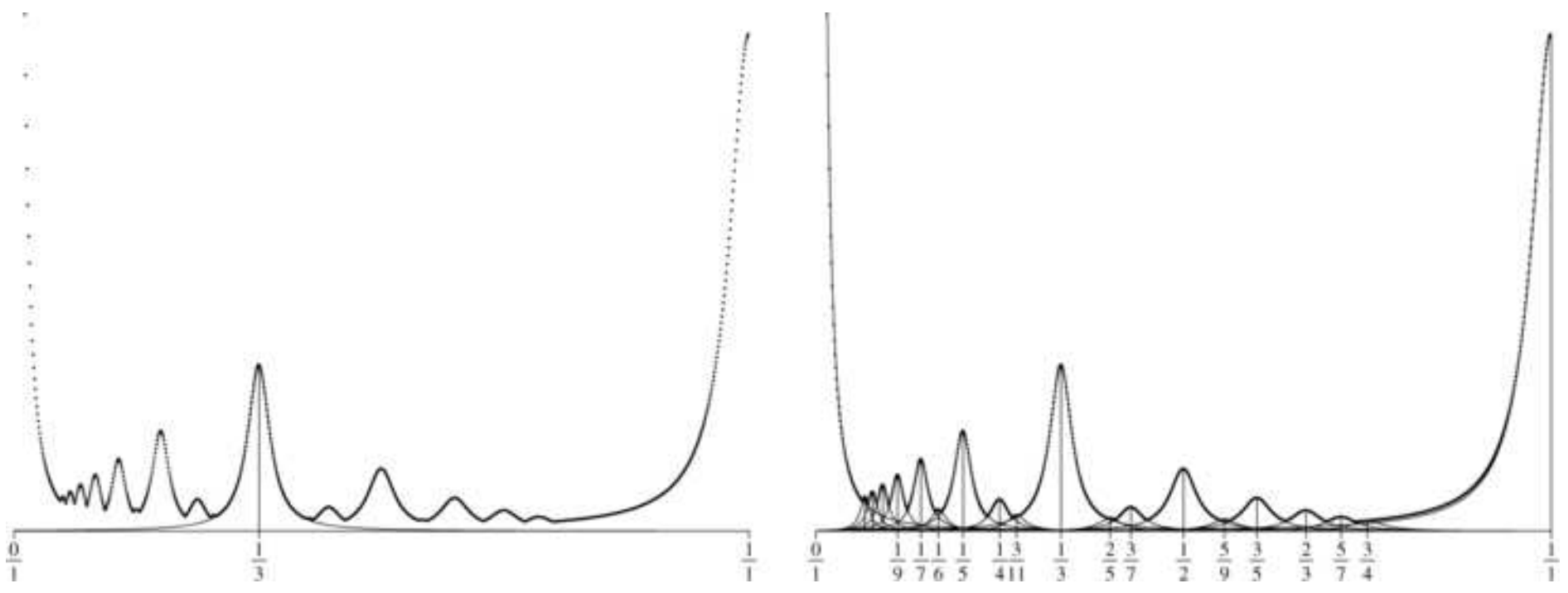

FIG. 3. *These plots show the same data from Figure 2 for the inverse spectral gap when $n = 1000$. On the left a specific bell curve is plotted on top of the spike at $k/n = 1/3$, and on the right such bell curves are plotted for all the simple fractions, including $\frac{0}{1}$ and $\frac{1}{1}$.*

*then the spectral gap is*

$$\operatorname{gap}(n,k) = (1+o(1))\frac{\pi^2 pq}{2n^2}\left(1 + \frac{q}{p}\frac{|k-(p/q)n|^2}{n}\right)$$
$$\times \begin{cases} 1, & \textit{if } p \equiv q \bmod 2, \\ 4, & \textit{if } p \not\equiv q \bmod 2. \end{cases}$$

*When $p/q = 0$, the spectral gap is instead $(1+o(1))2\pi^2(k+k^2)/n^3$, provided $k \le (n/2)^{2/3}$.*

Therefore, for each rational $p/q \neq 0$, the spectral gap ranges from order $1/n^2$ to order $1/n^{3/2}$ when $k/n$ is in a $n^{-1/4}$-neighborhood of $p/q$.

For any given $k$ and $n$, the spectral gap may be determined by the rational $k/n$ or perhaps by a nearby rational with smaller numerator and denominator. This leaves the question, for a given $k$ and $n$, assuming $n$ is large, how can we estimate the spectral gap? As Figure 3 indicates, the inverse spectral gap is determined by the rational whose bell is largest at $k/n$.

Let us define

$$\gamma(n,k,m) = \frac{\pi^2}{2n^2}\left(\frac{m^2 k}{n} + 4n\left\|\frac{m(n-k)}{2n}\right\|^2\right),$$
$$\gamma(n,k) = \min_{m\in\mathbb{Z}, m\neq 0} \gamma(n,k,m), \tag{1}$$

where $\|x\|$ is the number-theoretic norm $\|x\| = \operatorname{dist}(x, \mathbb{Z})$.

Our main theorem on the spectral gap, from which the previous theorems can be derived as corollaries, is the following:

THEOREM 4. *For all large enough $n$, for any $k$ such that $1 \le k \le n-1$, the spectral gap is*

$$\mathrm{gap}(n,k) = (1+o(1))\gamma(n,k),$$

*where the $o(1)$ goes to zero uniformly for the full range of $k$, when $n$ is large.*

This is a curious expression: there are terms in the min that do not correspond to eigenvalues, and there are eigenvalues that are not associated with terms in the min, but the sufficiently small terms do correspond to the eigenvalues with modulus sufficiently close to 1. Because the convergence holds uniformly for the full range of $k$, we can derive the previous theorems simply by studying $\gamma(n,k)$. For example, we have the following:

PROOF OF THEOREM 1. The min is achieved at $m = q$ when $p \equiv q \bmod 2$, and $m = 2q$ when $p \not\equiv q \bmod 2$. $\square$

An interesting consequence of Theorem 4 is the following:

THEOREM 5. *For $\alpha \in (0,1)$, we have*

$$\liminf_{n\to\infty} \gamma(n, \lfloor \alpha n \rceil) n^{3/2} \le 2\pi^2 \sqrt{\alpha}/\sqrt{5},$$

*and this bound is tight for a dense set of $\alpha$'s.*

PROOF. Any irrational number $x$ has infinitely many rational approximations $p/q$ for which $|x - p/q| < 1/(\sqrt{5}q^2)$, that is, $\|qx\| < 1/(\sqrt{5}q)$ (see Hardy and Wright [5], Section 11.8). Even for rational $x$, there are still infinitely many integers $q$ for which $\|qx\| < 1/(\sqrt{5}q)$. Let us set $x = (1-\alpha)/2$ and consider such an integer $q$ for which $q^2 \ge 1/\alpha$. We have, for $n \ge q^2$,

$$\begin{aligned}\gamma(n, \lfloor \alpha n \rceil, q) &= \frac{\pi^2}{2n^2}(q^2(\alpha + O(1/n)) + 4n\|q[(1-\alpha)/2] + O(q/n)\|^2) \\ &\le \frac{\pi^2}{2n^{3/2}}\left(\frac{q^2\alpha}{n^{1/2}} + \frac{4n^{1/2}}{5q^2} + \frac{O(1)}{n^{1/2}}\right).\end{aligned}$$

The optimal choice of $n$ is $n = \lceil 5q^4\alpha/4 \rceil > q^2$, for which we get

$$\gamma(n, \lfloor \alpha n \rceil, q) n^{3/2} \le \frac{\pi^2}{2}\left(\frac{4\sqrt{\alpha}}{\sqrt{5}} + o(1)\right),$$

where the $o(1)$ tends to 0 as $q \to \infty$, whence

$$\liminf_{n\to\infty} \gamma(n, \lfloor \alpha n \rceil) n^{3/2} \le \frac{2\pi^2\sqrt{\alpha}}{\sqrt{5}}.$$

There is a dense set of $\alpha$'s with continued fraction expansion that is eventually all 1's—these $\alpha$'s can be obtained from the golden ratio $\phi$ by the Möbius transformation $(a\phi+b)/(c\phi+d)$, where $a,b,c,d\in\mathbb{Z}$ and $ad-bc=1$ ([5], Section 10.11). For these $\alpha$'s, for any $\delta>0$, there are only finitely many rational approximations $p/q$ for which $\|qx\|<(1-\delta)/(\sqrt{5}q)$ ([5], Section 11.8), and so in the above calculations "$\leq$" may be replaced with $=(1+o(1))\times$, and so the above bound is tight. □

As the proof of the preceding theorem indicates, for a given $k$ and $n$, the optimal choice of $m$ for $\gamma(n,k,m)$ will be the denominator of a rational approximation $p/q$ to $(1-k/n)/2$. If $p/q$ is not a continued fraction convergent, then there is a continued fraction convergent $p_j/q_j$ which both more closely approximates $(1-k/n)/2$ and has smaller denominator $q_j<q$, so $\gamma(n,k,q_j)<\gamma(n,k,q)$. The best $q_j$ will never be larger than $O(n^{1/4})$, but it could be much smaller if $k/n$ has a particularly good rational approximation.

Remark. Diaconis and Saloff-Coste [2, 3] analyzed the mixing times of two related shuffles. Goel [4] and Jonasson [7] analyzed a variant in which one of the bottom $k$ cards is placed at the top of the deck.

Roadmap. In Section 2 we review the characteristic polynomial for this Markov chain and make some preliminary observations about the approximate locations of the eigenvalues. There are many eigenvalues that are almost on the unit circle. Then in Section 3 we prove Theorem 4 characterizing the spectral gap. Using Theorem 4, we calculate in Section 4 the spectral gap for extremal and for typical values of $k$, proving Theorem 2. In Section 5 we deduce Theorem 3 for the case when $k/n$ is well approximated by a simple rational.

**2. Eigenvalues for a single card.** The Markov chain for a single card is as follows: A card in position $i\in[1,n-k-1]$ moves to position $i+1$. The card in position $n-k$ moves to one of $\{1,n-k+1\}$, each with probability $\frac{1}{2}$. A card in position $i\in\{n-k+1,n-1\}$ moves to one of $\{i,i+1\}$, and the card in position $n$ moves to one of $\{1,n\}$, again each with probability $\frac{1}{2}$. Since most of the transition matrix is zero, it is easy to verify that its characteristic polynomial is

$$\lambda^{n-k}(\lambda-1/2)^k-(\lambda-1/2)^k/2-2^{-k-1},$$

and therefore, its eigenvalues satisfy

$$2\lambda^{n-k}-1=(2\lambda-1)^{-k}. \tag{2}$$

This equation for the eigenvalues first appeared in Jonasson's analysis [7], Section 4.1.2.

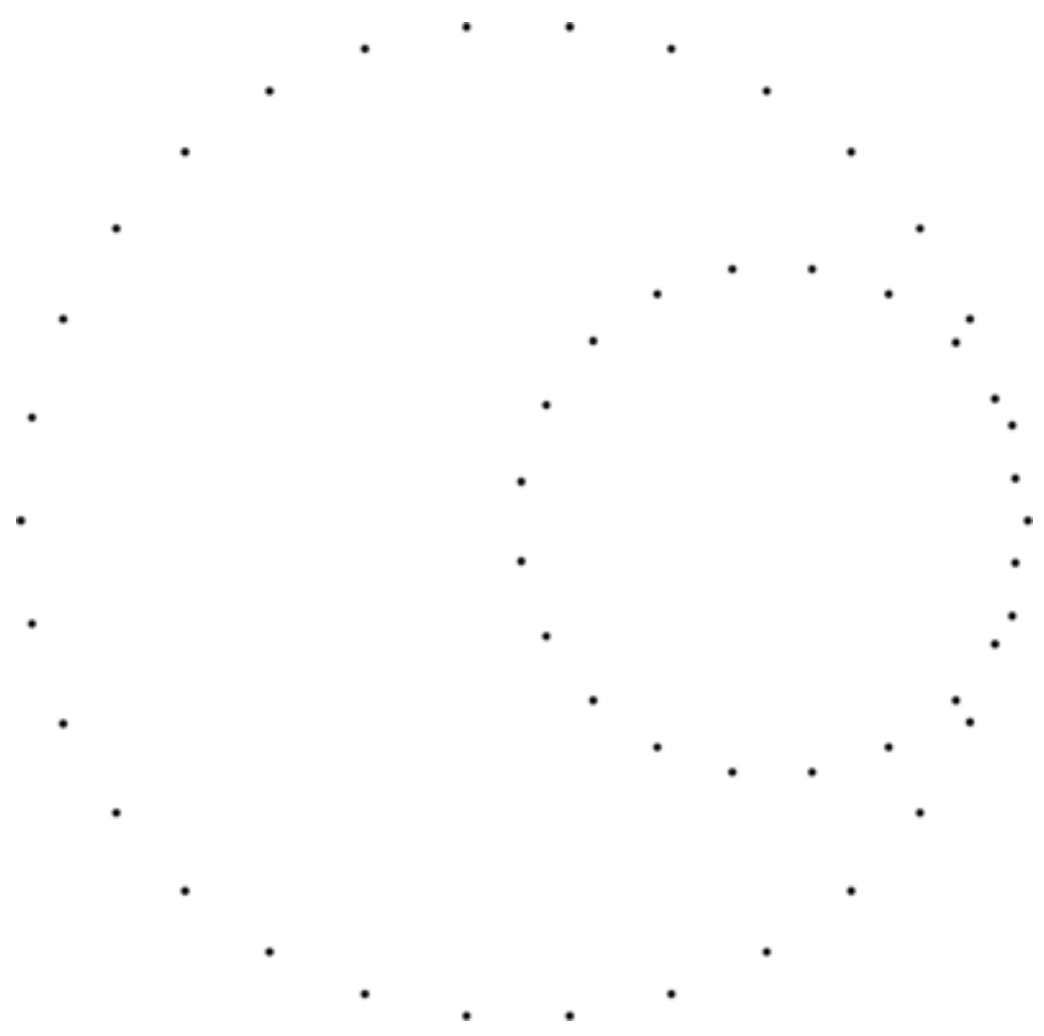

Fig. 4. *Eigenvalues of the* 50–20 *shuffle. In general, the eigenvalues will be close to the unit circle near the* $(n-k)$*th roots of unity, or near the circle of radius* $1/2$ *centered at* $1/2$*. Only eigenvalues close to* 1 *will be much closer to the unit circle than distance* $1/n$*.*

Proposition 6. *The spectral gap for the position of a single card is positive.*

Proof. Note that if $|z| \geq 1$, then $|2z-1| \geq 1$, with equality iff $z = 1$. Thus, if $|\lambda| \geq 1$, then the LHS of (2) is outside the unit circle and the RHS inside, and equality is possible iff $\lambda = 1$. Since the derivative of $(2\lambda^{n-k} - 1)(2\lambda - 1)^k - 1$ at $\lambda = 1$ is $2n$, the $\lambda = 1$ eigenvalue has multiplicity 1. □

There are actually two reasons that we call the shuffle the overlapping-cycles shuffle. The first is that the cards are shuffled according to one of the two cyclic permutations $(n, n-1, \ldots, 1)$ and $(n-k, n-k-1, \ldots, 1)$. The other reason is that the eigenvalues for the position of a single card form two overlapping cycles, as can be seen from Figure 4. To see this, note that if $\lambda$ is sufficiently far inside the unit circle, then $\lambda^{n-k} \approx 0$ and the LHS of (2) is approximately $-1$, so that

$$\lambda \approx \frac{(-1)^{1/k} + 1}{2},$$

where $(-1)^{1/k}$ is any $k$th root of $-1$ (sufficiently far from 1). If $|\lambda| \approx 1$ but $\lambda$ is sufficiently far from 1, then $(2\lambda - 1)^{-k} \approx 0$, so that (2) gives

$$\lambda \approx (1/2)^{1/(n-k)} (1)^{1/(n-k)},$$

where $(1)^{1/(n-k)}$ is any $(n-k)$th root of 1 (sufficiently far from 1). The most interesting range is where $\lambda \approx 1$, that is, where the two cycles overlap, and as we shall see in the next section, it is these eigenvalues that will ultimately determine the spectral gap of the Markov chain for a single card.

**3. Spectral gap for a single card.** We proceed to analyze the eigenvalues using equation (2) when $n$ is large. For an eigenvalue $\lambda$ of the overlapping cycles shuffle with parameters $n$ and $k$, we shall define $\varepsilon$, $a$ and $b$ to satisfy

$$\begin{aligned} \lambda &= \exp[-\varepsilon + i\pi a] && \text{where } -1 < a \le 1, \quad \text{and} \\ \lambda^{n-k} &= \exp[-(n-k)\varepsilon + i\pi b] && \text{where } -1 < b \le 1. \end{aligned} \tag{3}$$

At this point we restrict our attention to eigenvalues for which $\varepsilon = o(1/n)$. We shall see that there are always eigenvalues for which $\varepsilon = O(n^{-3/2})$, and for some values of $n$ and $k$ there are eigenvalues with much smaller $\varepsilon$. We start with a lemma showing that such eigenvalues satisfy $\lambda \approx 1$ and $\lambda^{n-k} \approx 1$, and derive successively more precise estimates of $\lambda$.

LEMMA 7. *With the notation (3), for any eigenvalue $\lambda$ for which $\varepsilon = o(1/n)$, we have $a = o(1)$ and $b = o(1)$.*

PROOF. We have $|2\lambda - 1| \ge 1 - 2\varepsilon + c(a)$, where $c(a) > 0$ for $a \ne 0$ is a positive constant depending on $a$. Then

$$|2\lambda^{n-k}| - 1 \le |2\lambda^{n-k} - 1| = |2\lambda - 1|^{-k} \le (1 - 2\varepsilon + c(a))^{-k}.$$

When $\varepsilon = o(1)$ and $a = \Theta(1)$, it follows that $|\lambda^{n-k}| < 1 - \tilde{c}(a)$, where $\tilde{c}(a)$ is a positive constant depending on $a$, which contradicts $\varepsilon = o(1/n)$.

We also have $|2\lambda^{n-k} - 1| > 1 - n\varepsilon + c(b)$ where $c(b) > 0$ is as above. Thus,

$$|2\lambda| - 1 \le |2\lambda - 1| = |2\lambda^{n-k} - 1|^{-1/k} < (1 - n\varepsilon + c(b))^{-1/k}.$$

When $\varepsilon = o(1/n)$ and $b = \Theta(1)$, it follows that $|\lambda| < 1 - \bar{c}(b)/k$, where $\bar{c}(b)$ is a positive constant depending on $b$, which contradicts $\varepsilon = o(1/n)$. $\square$

LEMMA 8. *With the notation (3), for any eigenvalue $\lambda$ for which $\varepsilon = o(1/n)$,*

$$\begin{aligned} |a\pi| &\le \sqrt{2n\varepsilon/k} + O((n\varepsilon/k)^{3/2}), \\ |b\pi| &\le \sqrt{2n\varepsilon} + O((n\varepsilon)^{3/2}). \end{aligned}$$

PROOF. Let $z = e^{-x+i\theta}$ where $x$ and $\theta$ are real and $x > 0$. By the triangle inequality, $\frac{1}{2} + |z - \frac{1}{2}| \ge |z| \ge 1 - x$, so

$$|2z - 1| \ge 1 - 2x. \tag{4}$$

But by the law of cosines,

$$
\begin{aligned}
|2z-1|^2 &= 1+4|z|^2-4|z|\cos\theta \\
&= 1+4-8x+O(x^2)-4(1-x+O(x^2))(1-\theta^2/2+O(\theta^4)) \\
&= 1-4x+2\theta^2+O(x^2)+O(x\theta^2)+O(\theta^4).
\end{aligned} \tag{5}
$$

From (2), using (4) with $z=\lambda$, we get

$$|2\lambda^{n-k}-1|^2 = |(2\lambda-1)^{-k}|^2 \le (1-2\varepsilon)^{-2k} = 1+4k\varepsilon+O(n^2\varepsilon^2), \tag{6}$$

and using (5) with $z=\lambda^{n-k}$, we get

$$|2\lambda^{n-k}-1|^2 = 1-4(n-k)\varepsilon+O(n^2\varepsilon^2)+2b^2\pi^2(1+O(n\varepsilon)+O(b^2)), \tag{7}$$

and upon combining (6) with (7), we get

$$4n\varepsilon+O(n^2\varepsilon^2) \ge 2b^2\pi^2(1+O(n\varepsilon)+O(b^2)).$$

Hence, taking square roots,

$$\sqrt{2n\varepsilon}\times(1+O(n\varepsilon)) \ge |b\pi|(1+O(n\varepsilon)+O(b^2)).$$

Since $\varepsilon=o(1/n)$, and $b^2=o(1)$ by Lemma 7, we learn $b=O(\sqrt{n\varepsilon})$, and hence,

$$|b\pi| \le \sqrt{2n\varepsilon}\times(1+O(n\varepsilon)).$$

From (2), using (4) with $z=\lambda^{n-k}$, we get

$$
\begin{aligned}
|2\lambda-1|^2 &= |2\lambda^{n-k}-1|^{-2/k} \\
&\le (1-2(n-k)\varepsilon)^{-2/k} \\
&= 1+\frac{4(n-k)\varepsilon}{k}+O((n\varepsilon/k)^2),
\end{aligned} \tag{8}
$$

and using (5) with $z=\lambda$, we get

$$|2\lambda-1|^2 = 1-4\varepsilon+O(\varepsilon^2)+2a^2\pi^2(1+O(\varepsilon)+O(a^2)), \tag{9}$$

and upon combining (8) with (9), we get

$$4n\varepsilon/k+O((n\varepsilon/k)^2)+O(\varepsilon^2) \ge 2a^2\pi^2(1+O(\varepsilon)+O(a^2)).$$

Hence, taking square roots,

$$\sqrt{2n\varepsilon/k}\times(1+O(n\varepsilon/k)) \ge |a\pi|(1+O(\varepsilon)+O(a^2)).$$

Since $\varepsilon=o(1/n)$, and $a^2=o(1)$ by Lemma 7, we learn $|a|=O(\sqrt{n\varepsilon/k})$, and hence,

$$|a\pi| \le \sqrt{2n\varepsilon/k}\times(1+O(n\varepsilon/k)). \qquad \square$$

LEMMA 9. *With the notation (3), for any eigenvalue $\lambda$ for which $\varepsilon = o(1/n)$, we can write*

$$\lambda = \exp\Big[\pi i m/n - \frac{\pi^2}{2n}(k(m/n)^2 + (m(1-k/n) \bmod 2)^2) + O(n^{1/2}\varepsilon^{3/2})\Big] \tag{10}$$

*for some integer $m = O(n^{3/2}\varepsilon^{1/2}/k^{1/2}) = o(n/k^{1/2})$, where mod 2 is interpreted so that $-1 < * \bmod 2 \le 1$.*

PROOF. We apply the following expansion for $z \approx 0$:

$$2e^z - 1 = \exp[2z - z^2 + O(z^3)]. \tag{11}$$

By Lemma 8, we know that eigenvalues $\lambda$ for which $\varepsilon = o(1/n)$ satisfy $|\log \lambda^{n-k}| = O(\sqrt{n\varepsilon})$, so

$$2\lambda^{n-k} - 1 = \exp[-2(n-k)\varepsilon + 2b\pi i + b^2\pi^2 + O((n\varepsilon)^{3/2})].$$

Lemma 8 also tells us $|\log \lambda| \le O(\sqrt{\varepsilon n/k})$, so

$$2\lambda - 1 = \exp[-2\varepsilon + 2a\pi i + a^2\pi^2 + O((\varepsilon n/k)^{3/2})],$$
$$(2\lambda - 1)^{-k} = \exp[2k\varepsilon - ka2\pi i - ka^2\pi^2 + O((\varepsilon n)^{3/2}/k^{1/2})].$$

For any eigenvalue $\lambda$, by (2), these expansions are equal:

$$-2(n-k)\varepsilon + 2b\pi i + b^2\pi^2 \equiv 2k\varepsilon - ka2\pi i - ka^2\pi^2 + O((n\varepsilon)^{3/2}) \bmod 2\pi i. \tag{12}$$

Comparing the real parts of (12), we obtain

$$-2n\varepsilon + b^2\pi^2 = -ka^2\pi^2 + O((n\varepsilon)^{3/2})$$
$$\varepsilon = \frac{\pi^2}{2n}(ka^2 + b^2) + O(n^{1/2}\varepsilon^{3/2}). \tag{13}$$

Comparing the imaginary parts of (12), we obtain

$$b \equiv -ka + O((n\varepsilon)^{3/2}) \bmod 1.$$

Of course, from (3) it is obvious that

$$b \equiv (n-k)a \bmod 2.$$

Subtracting these two relations, we obtain

$$0 \equiv na + O((n\varepsilon)^{3/2}) \bmod 1.$$

Since $na$ is very close to an integer, we define

$$m = \lfloor na \rceil = O(n^{3/2}\varepsilon^{1/2}/k^{1/2}) = o(n/k^{1/2}).$$

Then

$$\begin{aligned} a &= m/n + O(\varepsilon^{3/2}n^{1/2}),\\ b &= [m(1-k/n) + O(\varepsilon^{3/2}n^{3/2})] \bmod 2, \end{aligned}$$

where mod 2 is interpreted so that $-1 < * \bmod 2 \le 1$. We have then that

$$\begin{aligned} \lambda &= \exp[a\pi i - \varepsilon]\\ &= \exp\biggl[\pi i m/n - \frac{\pi^2}{2n}(k(m/n)^2 + (m(1-k/n) \bmod 2)^2) + O(n^{1/2}\varepsilon^{3/2})\biggr]. \end{aligned}$$

This estimate is valid for any eigenvalue $\lambda$ for which $\varepsilon = o(1/n)$. $\square$

LEMMA 10. *If $m$ is an integer such that $|m| = o(n/\sqrt{k})$ and $|m(1-k/n) \bmod 2| = o(1)$, then there is an eigenvalue of the form given in equation* (10).

PROOF. Let us define

$$\begin{aligned} b' &= m(1-k/n) \bmod 2 \qquad \text{with } -1 < b' \le 1,\\ X &= \frac{\pi^2}{2n}(k(m/n)^2 + b'^2),\\ \lambda(z) &= \exp[\pi i m/n - X + z],\\ f(z) &= (2\lambda(z)-1)^k(2\lambda(z)^{n-k}-1) - 1. \end{aligned}$$

Our aim is to show that there is some choice of $z \in \mathbb{C}$ such that $f(z) = 0$ and $|z| = O(n^{1/2}X^{3/2})$. Of course, if $f(z) = 0$, then $\lambda(z)$ is an eigenvalue. Our hypotheses imply $X = o(1/n)$, so $|z| = o(X)$, $\varepsilon = (1+o(1))X$, and $|z| = O(n^{1/2}\varepsilon^{3/2})$, so the existence of such a $z$ will prove the lemma.

Since

$$\lambda(z)^{n-k} = \exp[\pi i b' - (n-k)(X-z)],$$

the expansion (11) together with $|\pi b'| \le \sqrt{2nX} = o(1)$ imply that whenever $|z| = O(X)$ we have

$$2\lambda(z)^{n-k} - 1 = \exp[2\pi i b' - 2(n-k)(X-z) + \pi^2 b'^2 + O(n^{3/2}X^{3/2})].$$

The expansion (11) together with $|\pi m/n| \le \sqrt{2nX/k} = o(1)$ imply

$$2\lambda(z) - 1 = \exp[2\pi i m/n - 2(X-z) + \pi^2(m/n)^2 + O((Xn/k)^{3/2})].$$

Thus, modulo $2\pi i$,

$$\begin{aligned}\log(1+f(z)) &\equiv k[2\pi i m/n - 2(X-z) + \pi^2(m/n)^2 + O((Xn/k)^{3/2})] \\ &\quad + [2\pi i b' - 2(n-k)(X-z) + \pi^2 b'^2 + O(n^{3/2}X^{3/2})] \\ &= 2\pi i k m/n + \pi^2 k(m/n)^2 + 2\pi i b' - 2n(X-z) + \pi^2 b'^2 \\ &\quad + O(n^{3/2}X^{3/2}) \\ &= 2nz + 2\pi i k m/n + 2\pi i b' + O(n^{3/2}X^{3/2}) \\ &\equiv 2nz + O(n^{3/2}X^{3/2}) \qquad \pmod{2\pi i},\end{aligned}$$

whence

$$f(z) = 2nz + O(n^{3/2}X^{3/2}).$$

By considering $z$'s for which $|z|$ is a large enough constant times $n^{1/2}X^{3/2}$ [but still within the constraint $|z| = O(X)$], we may apply Rouché's theorem (see [1], Section 4.5.2) to deduce that there is a $z$ for which $|z| = O(n^{1/2}X^{3/2})$ and $f(z) = 0$. □

Recall the definition of $\gamma(n,k)$ from equation (1). We are now ready to prove Theorem 4, that is, that the spectral gap is $\text{gap}(n,k) = (1+o(1))\gamma(n,k)$. We actually show a bit more, that $\text{gap}(n,k) = \gamma(n,k) + O(n^{1/2}\gamma(n,k)^{3/2})$.

Proof of Theorem 4. Recall the definition

$$\gamma(n,k,m) := \frac{\pi^2}{2n}(k(m/n)^2 + (m(1-k/n) \bmod 2)^2),$$

which is the subexpression corresponding to $\varepsilon$ (ignoring the error terms) in equation (10).

For any real numbers $\alpha$ and $N \geq 1$, there is an integer $m$ so that $1 \leq m \leq N$ and $m\alpha$ is within $1/(\lfloor N \rfloor + 1) \leq 1/N$ of an integer (see [5], Section 3.8, Theorem 36). For our purposes, we are interested in $\alpha = \frac{1}{2}(1-k/n)$ and $N = 2^{1/2}n^{1/2}/k^{1/4}$. Then there is a positive integer $m \leq N$ for which

$$\gamma(n,k,m) \leq \frac{\pi^2}{2n}(k(N/n)^2 + (2/N)^2) = \frac{2\pi^2 k^{1/2}}{n^2} \leq 2\pi^2 n^{-3/2}.$$

Let $\tilde{m}$ denote the positive integer minimizing $\gamma(n,k,\tilde{m})$, that is, $\gamma(n,k,\tilde{m}) = \gamma(n,k)$. We have

$$\frac{\pi^2}{2n}k(\tilde{m}/n)^2 \leq \gamma(n,k,\tilde{m}) \leq \frac{2\pi^2 k^{1/2}}{n^2},$$

so $|\tilde{m}| \leq 2n^{1/2}k^{-1/4} = o(n/k^{1/2})$, and similarly,

$$\frac{\pi^2}{2n}(\tilde{m}(1-k/n) \bmod 2)^2 \leq \gamma(n,k,\tilde{m}) \leq \frac{2\pi^2 k^{1/2}}{n^2},$$

so $|\tilde{m}(1-k/n) \bmod 2| \le 2k^{1/4}n^{-1/2} = o(1)$. Since $\tilde{m}$ satisfies the hypotheses of Lemma 10, it gives us an eigenvalue $\lambda$ for which $\varepsilon = (1+o(1))\gamma(n,k,\tilde{m}) = (1+o(1))\gamma(n,k)$.

Since $\gamma(n,k) \le 2\pi^2 n^{-3/2}$, this proves Theorem 2 part (i).

Therefore, any eigenvalue $\lambda$ which achieves the spectral gap will satisfy $\varepsilon = o(1/n)$, so we can use Lemma 9 to express it according to equation (10) for some integer $m = m^*$. The spectral gap is then $(1+o(1))\gamma(n,k,m^*)$. We claim further that $m^* \ne 0$, since otherwise $\lambda = \exp[O(n^{1/2}\varepsilon^{3/2})]$, so $\varepsilon = O(n^{1/2}\varepsilon^{3/2}) = o(\varepsilon)$, so $\varepsilon = 0$ and $\lambda = 1$, which by Proposition 6 is not the eigenvalue achieving the spectral gap. Since there is an eigenvalue for which $\varepsilon = (1+o(1))\gamma(n,k) \le (1+o(1))\gamma(n,k,m^*)$, we conclude that the spectral gap is in fact $(1+o(1))\gamma(n,k)$. □

## 4. Spectral gap for extremal and typical values of $k$.

Recall that the upper bound on the spectral gap [Theorem 2 part (i)] was shown in the proof of Theorem 4. We proceed to prove the lower bound.

Proof of Theorem 2 part (ii). If for some nonzero $k$ and $m$ we have $\gamma(n,k,m) \le 4\pi^2 n^{-3}$, then

$$\frac{\pi^2}{2n} k(m/n)^2 \le \gamma(n,k,m) \le 4\pi^2 n^{-3},$$

and so $km^2 \le 8$. Similarly,

$$\frac{\pi^2}{2n}(m(1-k/n) \bmod 2)^2 \le \gamma(n,k,m) \le 4\pi^2 n^{-3},$$

so $|m(1-k/n) \bmod 2| \le \sqrt{8}/n$. Since $k \le 8$, for $n \ge 11$ this implies that $|m| \ne 1$, and hence, $m = \pm 2$ and $k \le 2$. Then for $n \ge 11$,

$$\gamma(n,k,m) = \frac{\pi^2}{2n}(k(2/n)^2 + (2k/n)^2) = \frac{2\pi^2}{n^3}(k+k^2),$$

and our assumption implies $k = 1$ and $\gamma(n,k,m) = 4\pi^2 n^{-3}$. Thus, whenever $n \ge 11$ we have $\gamma(n,k) \ge 4\pi^2 n^{-3}$, and so the result follows from Theorem 4. □

We now prove the last part of Theorem 2 which pertains to "typical" values of $k$.

Proof of Theorem 2 part (iii). We will show that for any $n \ge 1$ and any $\delta > 0$, for all but $4\delta^{2/3}n$ values of $k$, we have $\gamma(n,k) \ge (\pi^2/2)\delta n^{-3/2}$.

Let us start by assuming $k \ge \delta^{2/3}n$. If $m \ge \delta^{1/6}n^{1/4}$, then

$$\gamma(n,k,m) \ge \frac{\pi^2 m^2 k}{2n^3} \ge \frac{\pi^2 \delta}{2n^{3/2}}.$$

We also have

$$\gamma(n,k,m) \geq \frac{\pi^2}{2n}(m(1-k/n) \bmod 2)^2 = \frac{\pi^2}{2n^3}(mk^* \bmod 2n)^2,$$

where $k^* = n-k$ and mod denotes the centered mod. If $k^*$ is a uniformly random integer from the range $1 \leq k^* \leq 2n$, then $mk^* \bmod 2n$ is uniformly distributed among the multiples of $\gcd(m,2n)$. For any $x$ with $0 < x \leq n$, there are

$$\left(2\left\lceil \frac{x}{\gcd(m,2n)} \right\rceil - 1\right)\gcd(m,2n) \leq 2x + \gcd(m,2n)$$

values of $k^*$ $(1 \leq k^* \leq 2n)$ for which $|mk^* \bmod 2n| < x$, and a fortiori at most $x + \frac{1}{2}\gcd(m,2n)$ values of $k^*$ $(1 \leq k^* \leq n-1)$ for which $|mk^* \bmod 2n| < x$. Upon taking $x = \delta^{1/2}n^{3/4}$, and summing over $1 \leq m < \delta^{1/6}n^{1/4}$, we find that there are at most

$$\sum_{m=1}^{\lfloor \delta^{1/6}n^{1/4} \rfloor}\left(x + \frac{1}{2}\gcd(m,2n)\right) \leq \delta^{2/3}n + \frac{(\delta^{1/6}n^{1/4})(\delta^{1/6}n^{1/4}+1)}{4}$$

values of $k \geq \delta^{2/3}n$ and $m < \delta^{1/6}n^{1/4}$ for which $\gamma(n,k,m) < (\pi^2/2)\delta n^{-3/2}$. Thus, there are at most

$$2\delta^{2/3}n + \tfrac{1}{4}\delta^{1/3}n^{1/2} + \tfrac{1}{4}\delta^{1/6}n^{1/4} \tag{14}$$

values of $k$ for which $\gamma(n,k) < (\pi^2/2)\delta n^{-3/2}$.

If $\delta^{2/3}n \leq \frac{1}{4}$, then the bound (14) is less than one, so there are $0 \leq 4\delta^{2/3}n$ values of $k$ for which $\gamma(n,k) < (\pi^2/2)\delta n^{-3/2}$. If $\delta^{2/3}n \geq \frac{1}{4}$, then the bound (14) is $< 4\delta^{2/3}n$. □

**5. The narrow bells near simple rational $k/n$.** Let us suppose that $k/n$ is well approximated by a simple rational $p/q$ where $0 \leq p/q \leq 1$ and $\gcd(p,q) = 1$. Our aim here is to prove Theorem 3, which assumes that $p$ and $q$ are held fixed while $n \to \infty$.

Proof of Theorem 3. Let $c$ be defined by $k = (p/q)n + cn^{3/4}$. Recall

$$\gamma(n,k,m) = \frac{\pi^2}{2n^2}(m^2k/n + n(m(1-k/n) \bmod 2)^2).$$

Consider the term

$$n(m(1-k/n) \bmod 2)^2 = n(m(1-p/q-cn^{-1/4}) \bmod 2)^2.$$

As we saw in the proof of Theorem 1, if $c = 0$, this term is minimized at $m = Aq$, where

$$A = \begin{cases} 1, & \text{if } p \equiv q \bmod 2, \\ 2, & \text{if } p \not\equiv q \bmod 2. \end{cases}$$

When $m$ is divisible by $Aq$ and $mcn^{-1/4} \le 1$ we have

$$\gamma(n,k,m) = \frac{\pi^2}{2n^2} m^2 (k/n + c^2 n^{1/2}). \tag{15}$$

If $m$ is not divisible by $Aq$, then $mcn^{-1/4}$ would have to be as large as $(1-o(1))2/(Aq)$ before the second term of $\gamma(n,k,m)$ would be smaller than $O(1/n)$. But for $m$ this large, we have

$$\gamma(n,k,m) \ge \frac{\pi^2}{2n^2} \frac{k}{n} \frac{4n^{1/2}}{A^2q^2c^2} \tag{16}$$

on account of the first term.

If $c$ is small enough for the first term of (15) to be nonnegligible, then (16) is order $1/n$ and thus irrelevant to computation of $\gamma(n,k)$. If the first term of (15) is negligible, then we may compare (15) to (16), and we find that (15) determines $\gamma(n,k)$, provided

$$A^2q^2c^2 \le \frac{k}{n}\frac{4}{A^2q^2c^2},$$

$$c \le \frac{(4k/n)^{1/4}}{Aq}.$$

If $p \ne 0$ and $c \le (4p/q)^{1/4}/(Aq)$, we obtain

$$\gamma(n,k) = (1+o(1))\frac{\pi^2}{2n^2} A^2 q^2 (p/q + c^2 n^{1/2}).$$

When $p = 0$ we instead obtain

$$\gamma(n,k) = (1+o(1))\frac{2\pi^2}{n^3}(k+k^2)$$

for $k \le (n/2)^{2/3}$. □

## 6. Open problems.

1. We would like to understand the mixing time of the overlapping-cycles shuffle as a Markov chain on the symmetric group. Is it within a factor of $O(\log n)$ of the relaxation time of a single card?
2. For all large $n$ there is a $k$ of the form $k = n - (2/3^{1/4} + o(1))n^{1/4}$ such that $\gamma(n,k) = (2\pi^2/\sqrt{3} + o(1))n^{-3/2}$. In fact, it appears that for each real $\alpha \in (0,1]$ there is a $k$ of the form $k = (\alpha + o(1))n$ for which $\gamma(n,k) = (2\pi^2/\sqrt{3} + o(1))k^{1/2}n^{-2}$. Furthermore, it appears that $\gamma(n,k)$ can be no larger than this, as is illustrated in Figure 5. Does the inequality $\mathrm{gap}(n,k) \le (2\pi^2/\sqrt{3} + o(1))k^{1/2}n^{-2}$ hold?

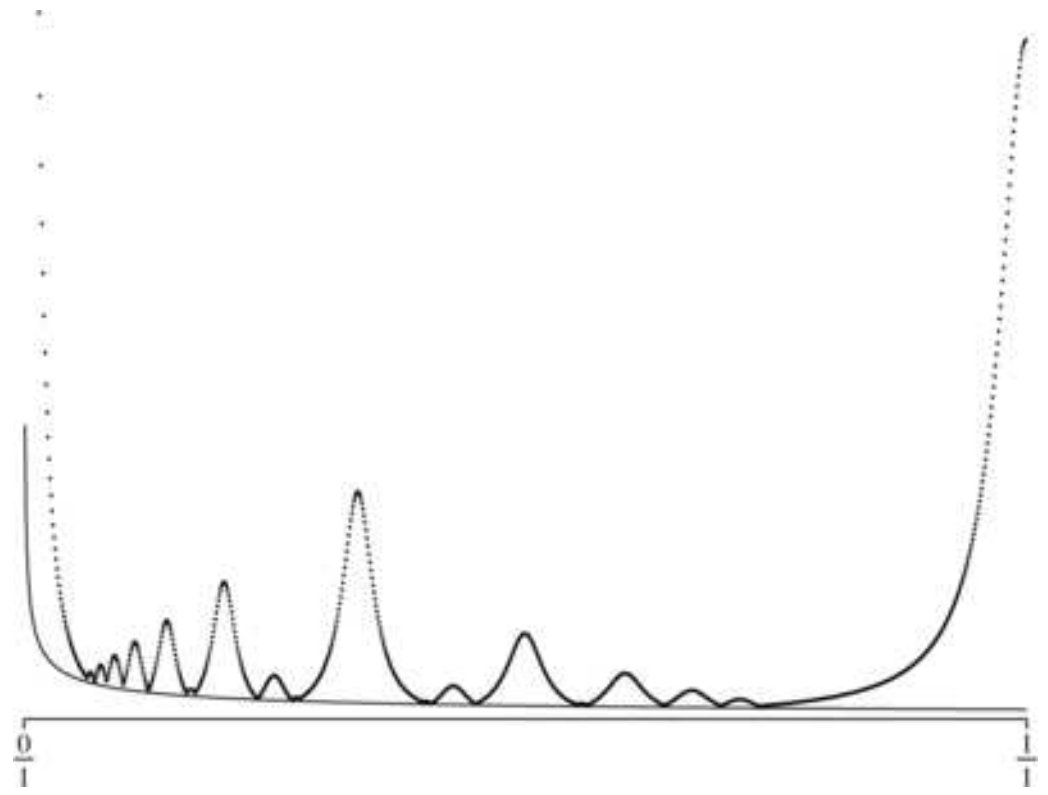

FIG. 5. *Shown here is the relaxation time for the position of the ♠A for the overlapping-cycles shuffle with parameters $n = 1000$ and $k$, together with a lower envelope on this function coming from a hypothetical bound $\gamma(n,k) \le (2\pi^2/\sqrt{3} + o(1))k^{1/2}n^{-2}$.*

**Acknowledgments.** This work originated in a walk by the beach during the probability summer school held at the University of British Columbia in June 2005. We are grateful to Ander Holroyd for useful discussions during that walk, and to Elizabeth Wilmer for help with simulations.

O. Angel
Department of Mathematics
University of Toronto
Toronto, Ontario
Canada M5S 2E4
E-mail: angel@math.toronto.edu
URL: http://www.math.toronto.edu/~angel

Y. Peres
D. B. Wilson
Microsoft Research
One Microsoft Way
Redmond, Washington 98052
USA
E-mail: peres@microsoft.com
dbwilson@microsoft.com
URL: http://research.microsoft.com/~peres
http://dbwilson.com